\documentclass[marginparwidth=30pt,a4paper,12pt]{article}

\usepackage[utf8]{inputenc}
\usepackage{graphicx}
\usepackage{amsmath}
\usepackage{gensymb}
\usepackage{multirow}
\usepackage{cleveref}
\usepackage{caption}
\usepackage{float}
\usepackage{xcolor}
\usepackage{authblk}

\usepackage{biblatex} %Imports biblatex package
\usepackage{comment}

\title{A discretization-convergent Level-Set-DEM}
\author[1]{Shai Feldfogel}
\author[2]{Konstantinos Karapiperis}
\author[3]{Jose Andrade}
\author[1]{David S. Kammer}
\affil[1]{Institute for Building Materials, ETH, Zurich, Switzerland}
\affil[2]{Department of Mechanical and Process Engineering, ETH, Zurich, Switzerland}
\affil[3]{Department of Mechanical and Civil Engineering, Caltech, Pasadena, California, USA}

\begin{document}

\begin{comment}

Why (purpose)
Convince anybody doing or that could do LS-DEM or DEM with surface discretization to use our approach

Who (audience): anybody doing LS-DEM, DEM community

What (Main messages)

0) Most important: All of you (audience): you should use our model!

1) 
\end{comment}

\maketitle

\begin{comment}
The ability of the level-set-DEM to seamlessly handle arbitrarily shaped grains and their contacts sets it apart from other DEM variants.
This ability is due to a discrete level-set representation of grains' volume and a node-based discretization of their bounding surfaces.
So far, the convergence properties of LS-DEM with refinement of these discretizations has not been studied.
Here, we show that the original LS-DEM diverges upon surface discretization refinement, and that this is due to its force-based discrete contact formulation.
We fix this by adopting a continuum-based contact formulation wherein the contact interactions are traction-based, and show that the adapted LS-DEM is fully discretization convergent.
Lastly, we evaluate the relative significance of the adapted LS-DEM formulation for different types of problems, and outline how it may be applicable in the broader context of DEM.
\end{comment}

\begin{abstract}
The recently developed level-set-DEM is able to seamlessly handle arbitrarily shaped grains and their contacts through a discrete level-set representation of grains' volume and a node-based discretization of their bounding surfaces.
Heretofore, the convergence properties of LS-DEM with refinement of these discretizations have not been examined.
Here, we examine these properties and show that the original LS-DEM diverges upon surface discretization refinement due to its force-based discrete contact formulation.
Next, we fix this issue by adopting a continuum-based contact formulation wherein the contact interactions are traction-based, and show that the adapted LS-DEM is fully discretization convergent.
Lastly, we discuss the significance of convergence in capturing the physical response, as well as a few other convergence-related topics of practical importance.

\begin{comment}
Lastly, we discuss the a few convergence-related topics of practical - the relation between surface-, and volume discretization refinement, the significance of convergence in capturing the physical response, and the importance of convergence in different types of problems commonly modeled with DEM.
\end{comment}

\end{abstract}

\section{Introduction} \label{sec:Introduction}
% What are TI structures, what is special about them, and why they are useful
% LS-DEM is a new DEM variant that overcomes the shape limitation of standard DEM. As such it a most promising computational approach for FATIS
Representing realistically-shaped grains has been a long-standing challenge in the Discrete Element Method (DEM) \cite{lai_fourier_2020,favier_shape_1999,peters_polyellipsoid_2009,hogue_shape_1998,feng_2d_2004}.
The Level-Set-Discrete-Element-Method (LS-DEM) was developed a few years ago to address this challenge \cite{lim_granular_2014,kawamoto_level_2016}, and it has since been applied to study various granular systems made up of realistically shaped grains.
Kawamoto et al. \cite{kawamoto_level_2016} simulated a triaxial compression test of Martian-like sand based on X-ray computed tomography of the actual grains used in the experiments.
Karapiperis et al. \cite{karapiperis_investigating_2020} investigated the constitutive elasto-plastic response and shape effects upon loading and unloading of Hostun sand.
Karapiperis et al. \cite{karapiperis_reduced_2020} studied the effects of reduced gravity on the strength of sand.
Harmon et al. \cite{harmon_level_2020} studied the effects of brittle particle breakage in the context of realistic particles going through a Jaw crusher, in an odeometric set-up, and in wall demolition.
Karapiperis et al. \cite{karapiperis_stress_2022} studied the formation and evolution of stress transmission paths in entangled granular assemblies, and their macroscopic mechanical properties.

LS-DEM's ability to represent arbitrarily shaped grains and resolve contact between such grains is due to two discretizations - a node-based discretization of grain boundary and a discretized Level-Set representation of grain volume.
The refinement of the surface discretization (SD) is defined by the distance between surface nodes and it determines the accuracy with which the contact regions are spatially resolved.
The refinement of the volume discretization (VD) is defined by the step size of the volume grid on which the Level-Set in calculated, and it determines the accuracy with which the magnitude and the direction of the penetration at a contact node are evaluated.

Convergence upon discretization refinement is important because it assures that results are discretization independent and quantitatively meaningful.
It is therefore a common requirement from discretization-based numerical models such as  the Finite Difference Methods and the Finite Element Method (FEM).
In this study, we address of the convergence properties of the discretization-based LS-DEM, which have not yet been addressed in the literature.
\begin{comment}
, likely because it grew out of DEM, where discretization and convergence are not relevant
This may be because LS-DEM grew out of the field of discrete elements, where discretization, and hence convergence, are not relevant.
Nevertheless, convergence is important for LS-DEM as it is in other methods as it assures that results are discretization independent and quantitatively meaningful.
\end{comment}

\begin{figure}[H] 
    \includegraphics[width=1\textwidth]{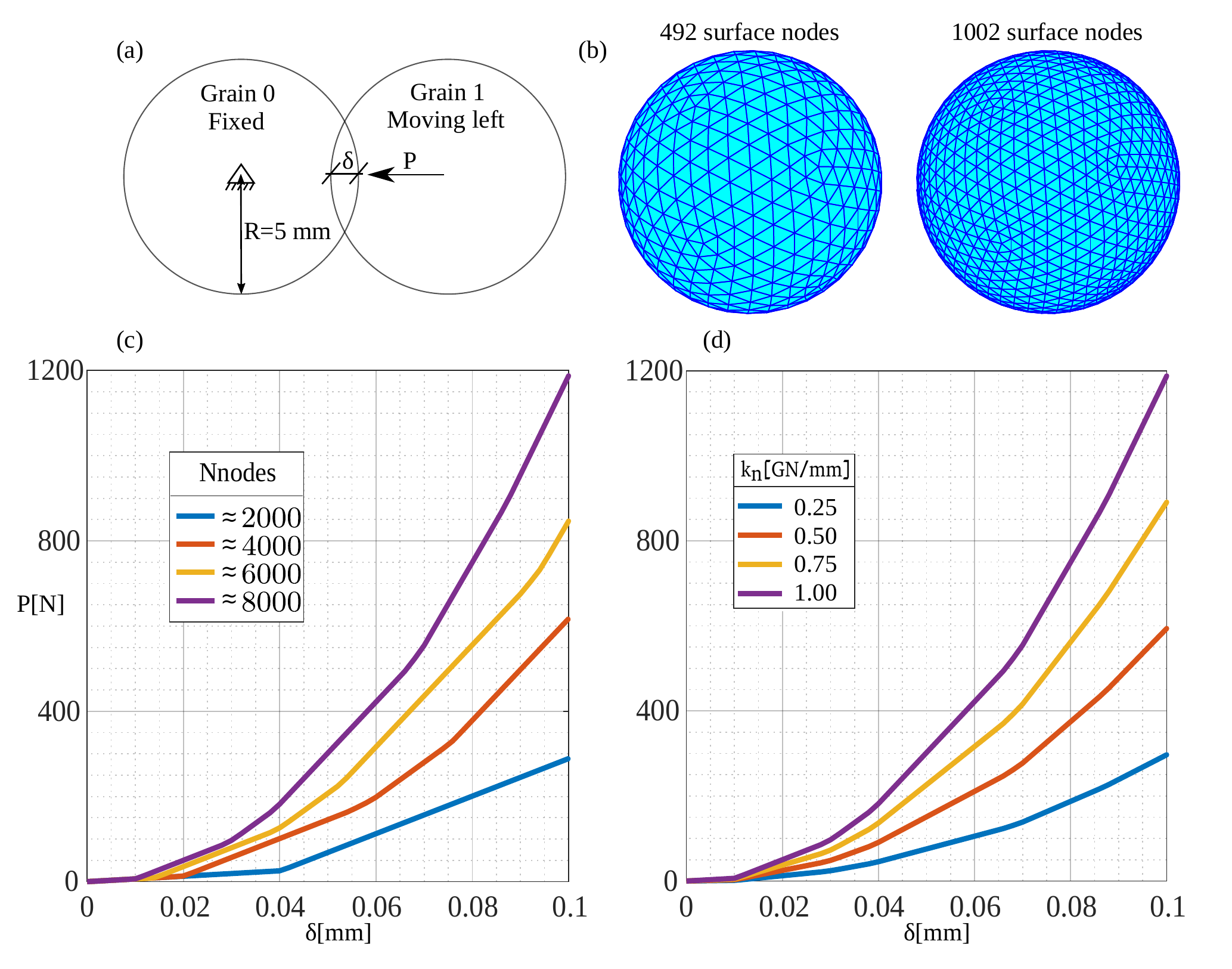} 
    \centering
    \caption{Two-sphere central-compression with original LS-DEM: (a) Configuration; (b) two surface discretizations; (c) P-$\delta$ curves diverge in the original LS-DEM (d) mesh refinement equivalent to increasing $k_n$.}
    \label{fig:twoSphereCompression_original}
\end{figure}

The convergence of LS-DEM upon SD refinement is considered in Fig \ref{fig:twoSphereCompression_original}.
Fig. \ref{fig:twoSphereCompression_original}(a) depicts a case of central-compression between two-spheres.
$P$ and $\delta$ denote the compression force and the penetration depth, respectively.
Fig. \ref{fig:twoSphereCompression_original}(b) shows two SDs differing in the total number of surface nodes (Nnodes).
As in all the examples in this study, the surface of the grains is represented by a triangle mesh, whose nodes constitute the SD.
Fig. \ref{fig:twoSphereCompression_original}(c) shows $P-\delta$ curves from four analyses with a different Nnodes.
Fig. \ref{fig:twoSphereCompression_original}(d) shows $P-\delta$ curves from analyses with Nnodes=8,000 and with four linearly-spaced  values of the penalty parameter $k_n$.
Fig. \ref{fig:twoSphereCompression_original}(c) shows that the response scales approximately linearly with SD refinement, which means that LS-DEM diverges upon SD refinement.
Fig. \ref{fig:twoSphereCompression_original}(d) shows that the response also scales approximately linearly with $k_n$, implying that SD refinement is equivalent to increasing $k_n$.

The divergent behavior observed in Fig. (c) underscores the need to adapt LS-DEM so that it becomes discretization convergent and it raises several questions of interest:
\begin{enumerate}
    \item Why does LS-DEM divergence upon SD refinement, and how can this divergence be fixed?
    \item How does the relation between the SD and VD refinements affect convergence?
    \item What is the physical significance of convergence in terms of correctly capturing the mechanical response?
    %\item How does convergence depend on the type of the problem considered?
\end{enumerate}
The main objectives of the present study are to fix LS-DEM's divergence issue and to address these questions.

Next, in section \ref{sec:Methodology} we identify the cause for LS-DEM divergence and propose an adapted contact formulation that eliminates this divergence.
In section \ref{sec:Results and discussion}, we examine the convergence properties of the adapted formulation upon SD refinement and address the aforementioned open questions.
We summarize the study in Section \ref{sec:Summary and outlook}.

\section{Methodology} \label{sec:Methodology}
The reason why the original LS-DEM diverges upon SD refinement is discussed in Sec. \ref{subsection:the root cause of divergence}; The elimination of this divergence is described in \ref{subsection:How we eliminate LSDEM's divergence}; The adapted formulation is detailed in \ref{subsection:Adapted LSDEM formulation}.

\subsection{Why does the original LS-DEM diverge?} \label{subsection:the root cause of divergence}

The root cause of LS-DEM's divergent behavior observed in Figure \ref{fig:twoSphereCompression_original}(a) is that the inter-grain penetration stiffness, which controls the global strength and stiffness, scales linearly with the number of surface nodes.

\begin{figure}[H]
    \includegraphics[width=1\textwidth]{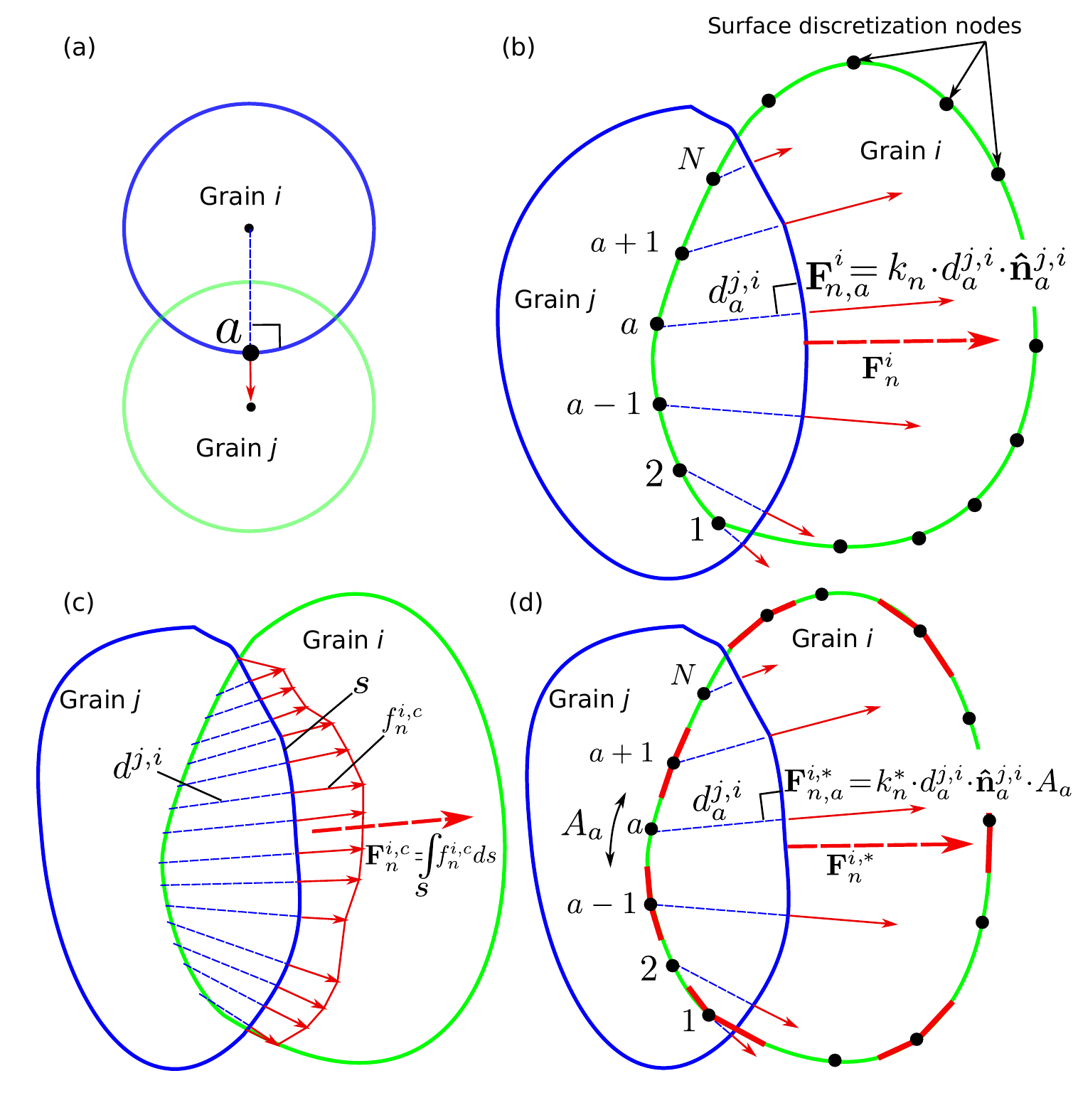} 
    \centering
    \caption{The adapted LS-DEM overcomes the divergent behavior of the original LS-DEM by adopting a continuum-based contact description: (a) ordinary DEM; (b) original LS-DEM; (c) continuum contact model; and (d) adapted LS-DEM (penetrations are grossly exaggerated for illustrative purposes)}
    \label{fig:methodology}
\end{figure}

To see this in more detail, consider Fig. \ref{fig:methodology}.
In ordinary DEM, the grains are spherical and each contact comprises a single penetration point $a$ as illustrated in Fig \ref{fig:methodology}(a).
In LS-DEM, due to the surface discretization, there are generally multiple contact nodes between two contacting grains, as shown in Figure \ref{fig:methodology}(b).
At each node, the normal force acting on grain $i$ by grain $j$ is, similarly to ordinary DEM:
\begin{equation}\label{eq:original LSDEM Fna}
    \mathbf{F}^{i}_{n,a} = k_n \cdot d^{j,i}_{a} \cdot \mathbf{\hat{n}}^{j,i}_{a} 
\end{equation}
where $k_n$ is the normal penalty parameter with dimensions [force/penetration],  \(d^{j,i}_a\) and \(\hat{\mathbf{n}}^{j,i}_a\) denote, respectively, the penetration depth and the unit normal to grain \textit{j} at penetration point \textit{a}, and where the summation is carried over the surface nodes on the boundary of grain $j$.

\begin{comment}
and are calculated, 
respectively, as the value and the gradient of the level-set of grain $i$
\end{comment} 
\begin{comment}
, commonly thought of as a discrete linear spring whose deformation corresponds to the penetration depth.
\end{comment}

The resultant normal force $\mathbf{F}^{i}_{n}$ acting on grain $i$ by grain $j$, marked in Fig. \ref{fig:methodology} by a dashed arrows, is the of vector sum of nodal forces over the \textit{N} contact nodes:
\begin{equation}\label{eq:original LSDEM Fn}
    \mathbf{F}^{i}_{n} = \sum_{a=1}^{N} \mathbf{F}^{i}_{n,a} = \sum_{a=1}^{N} k_n \cdot d^{j,i}_{a} \cdot \mathbf{\hat{n}}^{j,i}_{a} 
\end{equation}

Then, defining the average penetration $\overline{d}^{j,i}_{a}\equiv\frac{\sum_{a=1}^{N} d^{j,i}_{a}}{N}$ and the grain penetration stiffness $k^{grain}$ as $\frac{|\mathbf{F}^{i}_{n}|}{\overline{d}^{j,i}_{a}}$, it follows that $k^{grain}$ is linear in both $N$ and $k_n$:
\begin{equation}\label{eq:penetration_stiffness_original_LSDEM}
    k^{grain} \equiv \frac{|\mathbf{F}^{i}_{n}|}{\overline{d}^{j,i}_{a}} =
  \frac{\sum_{a=1}^{N} |k_n \cdot d^{j,i}_{a} \cdot \mathbf{\hat{n}}^{j,i}_{a}|}{\frac{\sum_{a=1}^{N} d^{j,i}_{a}}{N}} = N \cdot k_n \cdot \frac{\sum_{a=1}^{N} |d^{j,i}_{a} \cdot \mathbf{\hat{n}}^{j,i}_{a}|}{\sum_{a=1}^{N} d^{j,i}_{a}} = N \cdot k_n
\end{equation}
The distinction between the grain penetration stiffness $k^{grain}$ which refers to interaction between grains as a whole and the normal penalty parameter $k_n$ which refers to contact interactions at the node level does not exist in ordinary DEM. There, $N$=1 and $k^{grain}$ and $k_n$ are therefore one and the same.

The linearity of $k^{grain}$ with $k_n$ in Eq. (\ref{eq:penetration_stiffness_original_LSDEM}) is reflected in Fig.  \ref{fig:twoSphereCompression_original}(b), and it is common to all DEM variants.
It is useful because it allows to calibrate $k_n$ as a correlate of the elastic modulus.
In contrast, the linearity of $k^{grain}$ with $N$ is unique to LS-DEM and it is problematic because, as the SD is refined, $N$ and $k^{grain}$ increase proportionally, leading to the divergence observed in Fig. \ref{fig:twoSphereCompression_original}(a).

\subsection{How we eliminate LS-DEM's divergence} \label{subsection:How we eliminate LSDEM's divergence}
The key to fixing LS-DEM's divergence is to eliminate the  discretization dependence of the penetration stiffness expressed in Eq. (\ref{eq:penetration_stiffness_original_LSDEM}).
To do that, we adopt the continuum-based approach described in Figure \ref{fig:methodology}(c), where contact forces are integrals of surface tractions over the contact area $s$, rather than sums of nodal forces.
The dimensions of the continuum-based penalty parameter $k_n^{*}$ are those of an elastic foundation [traction/penetration], differently from  $k_n$'s [force/penetration]. 
This approach provides the resultant normal force \textbf{F}$^{i,c}_{n}$ and the commensurate penetration stiffness $k^{grain}$, with the necessary linkage to and bounding by the finite contact area $s$, bounding which the original formulation expressed in Equation (\ref{eq:original LSDEM Fn}) lacks.

The continuum-based resultant normal force $\mathbf{F}^{i,c}_n$ reads:
\begin{equation} \label{eq:analytical Fn}
    \mathbf{F}^{i,c}_n = \int_s d\mathbf{F}^{i,c}_n 
\end{equation}
where the superscript $^{,c}$ denotes the continuum-based approach, $s$ denotes the contact area between the grains, and the elemental $d\mathbf{F}^{i,c}_n$ reads:
\begin{equation} \label{eq:analytical dFn}
    d\mathbf{F}^{i,c}_n = \mathbf{f}^{i,c}_n ds = k^{*}_n \cdot d^{j,i} \cdot \mathbf{\hat{n}}^{j,i} ds
\end{equation}
where $\mathbf{f}^{i,c}_n$ are the normal tractions, \(d^{j,i}\) the normal penetrations, and \(\mathbf{\hat{n}}^{j,i}\)\ the surface normals.
\begin{comment}
$k^{*}_n$ represents an elastic foundation, with dimensions of [traction/penetration], differently from the discrete spring $k_n$ [force/penetration] in Equations (\ref{eq:original LSDEM Fna}-\ref{eq:penetration_stiffness_original_LSDEM}). 
\end{comment}

We evaluate the integral in Equation (\ref{eq:analytical Fn}) numerically as a sum of nodal contributions at contacting nodes, as illustrated in Fig. \ref{fig:methodology}(d).
The resultant normal force in the adapted LS-DEM $\mathbf{F}^{i,*}_n$ therefore reads:
\begin{equation} \label{eq:adapted LSDEM Fn}
    \mathbf{F}^{i,*}_n = \sum_{a=1}^{N} \mathbf{F}^{i,*}_{n,a}
\end{equation}
where the superscript $^{,*}$ denotes the adapted formulation.
The nodal normal force $\mathbf{F}^{i,*}_{n,a}$, which is the finite analog of $d\mathbf{F}^{i,c}_n$ from Equation (\ref{eq:analytical dFn}), reads: 
\begin{equation} \label{eq:adapted LSDEM Fna}
    \mathbf{F}^{i,*}_{n,a} = k^{*}_n \cdot d^{j,i}_a \cdot \mathbf{\hat{n}}^{j,i}_a \cdot A_a
\end{equation}
where $A_a$ is the tributary area of node $a$, see Fig. \ref{fig:methodology}(d).
\begin{comment}
As in Equation (\ref{eq:original LSDEM Fna}), $d^{j,i}$ and $\mathbf{\hat{n}}^{j,i}$ are evaluated from the level-set of grain $i$, and the summation is carried over the contacting surface nodes on the boundary of grain $j$.
\end{comment}

The penetration stiffness $k^{int,*}$ in the adapted model then reads:
\begin{equation}\label{eq:penetration_stiffness_adapted_LSDEM}
    k^{int,*} \equiv \sum_{a=1}^{N} \frac{|\mathbf{F}^{i,*}_{n,a}|}{d^{j,i}_{a}} =
              \sum_{a=1}^{N} \frac{|k^{*}_n \cdot d^{j,i}_{a} \cdot \mathbf{\hat{n}}^{j,i}_{a} \cdot A_a|}{d^{j,i}_{a}} = \sum_{a=1}^{N} k^{*}_n \cdot A_a = k^{*}_n \cdot N \cdot A_a
\end{equation}
The "corrective" factor $A_a$ and the elastic-foundation nature of $k_n^{*}$ are the two keys to resolving the divergence issue.
By taking \textit{A}$_a$ equal to the total surface area of the grain divided by the total number of surface nodes, an increase in $N$ with SD refinement is off-set ("corrected") by a proportional decrease in \textit{A}$_a$, so that the sensitivity of $k^{int,*}$ on the SD becomes minor, and the global response converges upon SD refinement, as will be shown in the next section. 
\begin{comment}
\footnote{The SD refinement should be above a certain minimum to reasonably capture the contact area}
\end{comment}

The tangential contact formulation requires a treatment similar to that of the normal contact formulation, for similar reasons.
Comparing Eqs. (\ref{eq:original LSDEM Fna}) and (\ref{eq:adapted LSDEM Fna}), the modification of the normal nodal force in the adapted formulation amounts to replacing $k_n$ with the $k^{*}_n \cdot A_a$.
The modification of the nodal tangential force in the adapted formulation is analogous to this replacement.
In the original LS-DEM (see Equation (15) in \cite{kawamoto_level_2016}), the elastic tangential force at node $a$, which is the basis for the commonly-used Coulomb friction law, is based on the product $k_s \cdot \Delta\mathbf{S}_a$, where $k_s$ [force/displacement] is the tangential stiffness analogous to $k_n$ and where $\Delta\mathbf{S}_a$ is the (incremental) tangential slip vector analogous to the penetration vector $d^{j,i}_{a} \cdot \mathbf{\hat{n}}^{j,i}_a$. 
Hence, the modification of the nodal shear force amounts to replacing $k_s$ with $k^{*}_s \cdot A_a$, where the distributed tangential stiffness $k^{*}_s$ has dimensions [traction/penetration]. 
\subsection{Adapted LS-DEM formulation} \label{subsection:Adapted LSDEM formulation}
With the modifications to the contact formulation described above, the adapted LS-DEM's formulation is neatly obtained from the original one by replacing $k_n$ and $k_s$ in Equations (5) and (9) from \cite{kawamoto_level_2016} with $k^{*}_n \cdot A_a$ and $k^{*}_s \cdot A_a$, respectively.
Regarding the computational aspect, it is pointed out that the additional computational cost in the adapted formulation is effectively zero. 
It consists merely of calculating $A_a$ once for each node in the pre-processing stage and appending it to $k_n^{*}$ and $k_s^{*}$ whenever the node is in contact.

\section{Results and discussion} \label{sec:Results and discussion}
The convergence properties of the adapted LS-DEM in different configurations are evaluated and discussed in sections \ref{subsection:Two-sphere central-compression}-\ref{subsection:A topologically interlocked slab}.
The other questions posed at the end of the Introduction are addressed in sections \ref{subsection:The relation between SD and VD refinement}-\ref{subsection:The physical significance of convergence}.
\subsection{Two-sphere central-compression} \label{subsection:Two-sphere central-compression}
\begin{figure}[H] 
    \includegraphics[width=1\textwidth]{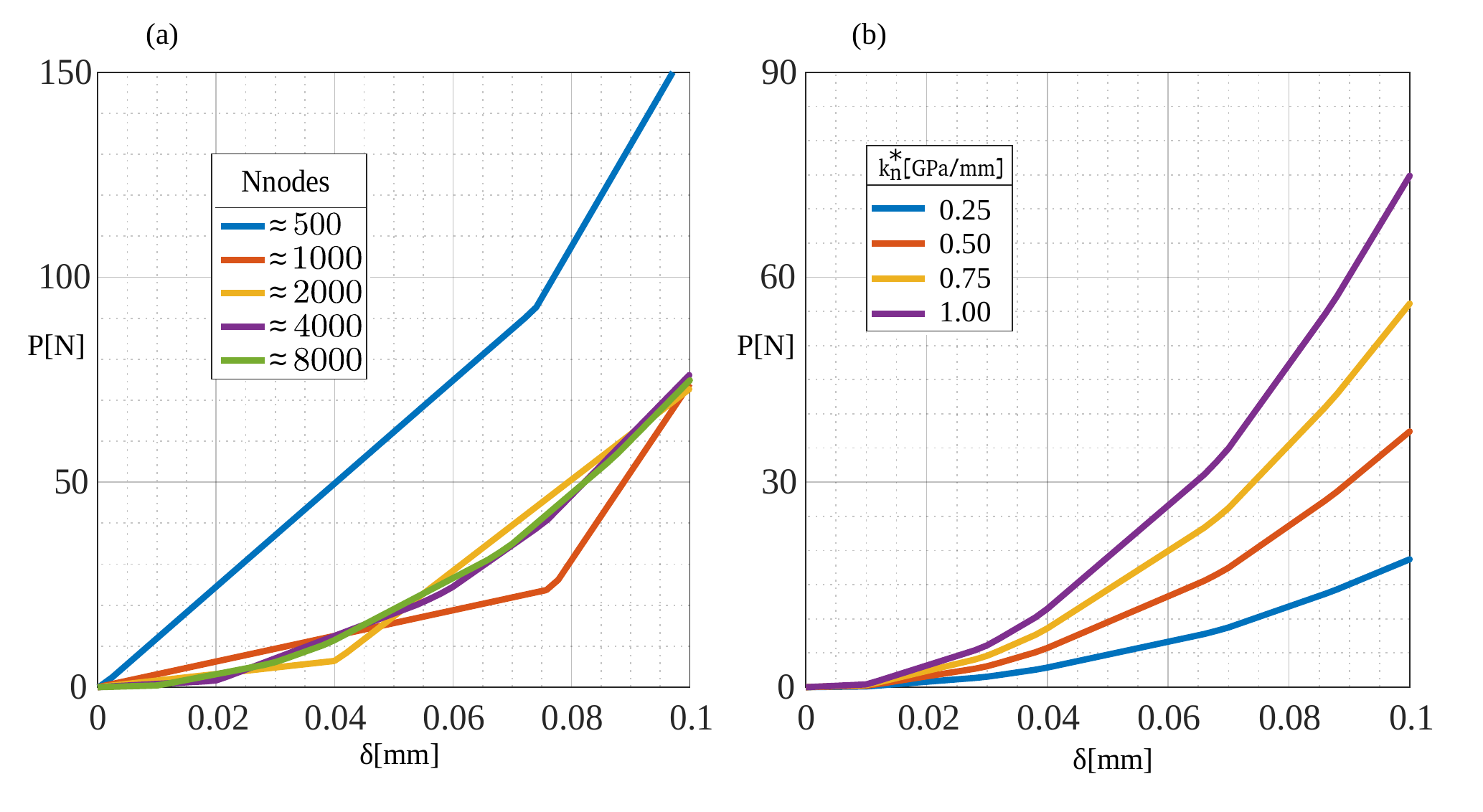} 
    \centering
    \caption{The two-sphere central-compression case with the adapted LS-DEM: (a) the P-$\delta$ curves converge upon SD refinement; (b) the response scales linearly with $k_n^{*}$ in accordance with Eq. (\ref{eq:penetration_stiffness_adapted_LSDEM}).}
    \label{fig:twoSphereCompression_adapted}
\end{figure}

We begin the evaluation of the adapted LS-DEM's convergence properties by revisiting the two-sphere configuration described in Fig. \ref{fig:twoSphereCompression_original}.
Fig. \ref{fig:twoSphereCompression_adapted}(a) depicts five $P-\delta$ curves obtained with the adapted LS-DEM and with the same SD's as those used for Fig. \ref{fig:twoSphereCompression_original}(c).
The $P-\delta$ curves in Fig. \ref{fig:twoSphereCompression_adapted}(a) show that, in contrast to the divergent behavior of the original LS-DEM seen in Fig. \ref{eq:penetration_stiffness_original_LSDEM}(c), the adapted formulation converges upon SD refinement. 
Fig. \ref{fig:twoSphereCompression_adapted}(b) shows the linearity of the response with $k_n^{*}$, reflecting the linearity of $k^{int,*}$ with $k_n^{*}$ in accordance with Eq. (\ref{eq:penetration_stiffness_adapted_LSDEM}).
The convergence of the adapted formulation and the linearity of the response with $k_n^{*}$ in accordance with Eq. (\ref{eq:penetration_stiffness_adapted_LSDEM}) show that cause of divergence has indeed been correctly identified and eliminated as described in Sec. \ref{sec:Methodology}.

%\subsection{A pill box} \label{subsection:A pill box}
%Next, we consider the assembly of spheres/pills described in Fig. 

\subsection{A topologically interlocked slab} \label{subsection:A topologically interlocked slab}
Moving on from the simple two-sphere case to a more realistic configuration with multiple particles, we consider next the topologically interlocked slabs studied experimentally in \cite{mirkhalaf_toughness_2019}.
Topologically interlocked slabs are made of un-bonded building blocks that hold together by contact and friction thanks to their interlocking shapes.
Their analysis with LS-DEM represents the latter's novel application to structural analysis problems, where convergence is particularly important. 

\begin{figure}[H] 
    \includegraphics[width=1\textwidth]{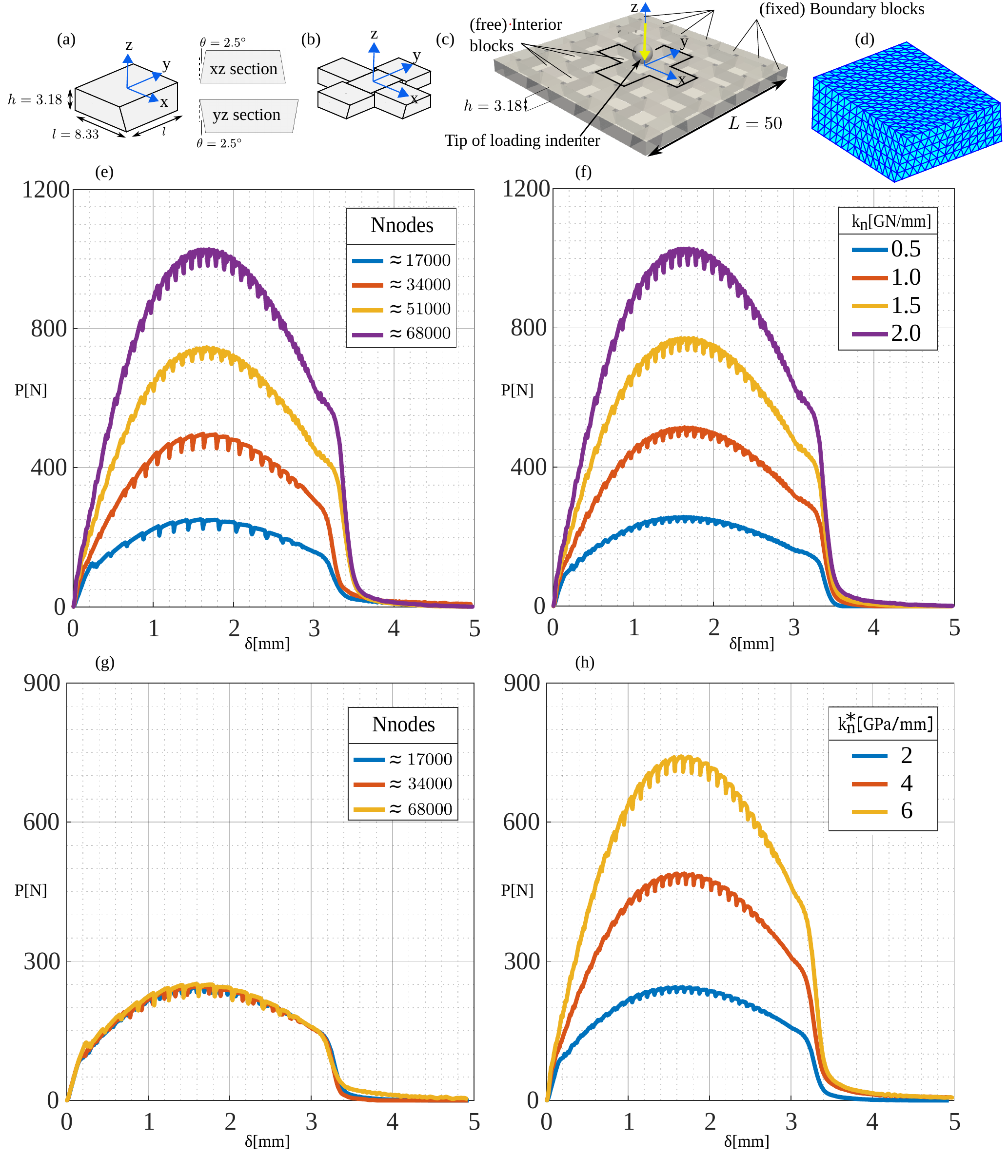} 
    \centering
    \caption{Topologically interlocked slab: (a-d) Configuration; (e-h) P-$\delta$ curves diverge in the original LS-DEM (e) and converge in the adapted LS-DEM (g) upon discretization refinement; the response scales linearly with the penalty parameters $k_n$ and $k_n^{*}$ (f,h).}
    \label{fig:TIS}
\end{figure}

Fig. \ref{fig:TIS}(a) shows a typical planar-faced block that makes up the interlocking slab and its cross-sections in the xz and yz planes.
The bottom face of the block is a square with side length $l=8.33$ mm, and the angle of inclination of its sloping lateral faces is $\theta$=2.5$^{o}$.
Fig. \ref{fig:TIS}(b) shows the basic 5-block interlocked cell formed by surrounding a block by four similar ones rotated with respect to it by 90$^{\circ}$ with respect to the z axis.
Fig. \ref{fig:TIS}(c) shows the entire slab with the contour of a basic cell around the central block marked in black.
The slabs' dimensions are 50 x 50 x 3.18 mm, and they consist of boundary blocks that are fixed by an external constraint along the edges and of internal blocks.
Fig. \ref{fig:TIS}(d) shows the triangular mesh constitutes the SD of the central block. 

The slab was loaded by a concentrated load $P$ in the negative z direction.
The load was affected by prescribing a constant velocity to a spherical loading indenter with a 2.5 mm radius.
$P$ and the corresponding displacement $\delta$ are indicated by a yellow arrow in the negative z direction in Figure \ref{fig:TIS}(c).
Prior to the indentation loading phase, the assembly was subjected to gravity until it reached a relaxed state, i.e., until the kinetic energy lowered to effectively zero.
The relaxed positions and rotations of the blocks were taken as the initial conditions for the concentrated load phase.
\begin{comment}
To expedite the analyses, the loading speed was taken as high as possible, but always low enough to ensure stability of the time stepping algorithm, and to avoid any inertial effects.
The loading rate values ranged between 3-6 mm/sec.
\end{comment}
The material density was taken equal to $2.5\cdot10^{-6}\frac{kg}{mm^3}$, and a friction coefficient $\mu$=0.23 was used, in accordance with \cite{mirkhalaf_toughness_2019}. 

\begin{comment}
The mesh's nodes are the seeded nodes and the LS values represent the surface discretization and was constructed with the seeded surface nodes as its vertices; and (3) the block's level-set function was calculated within a bounding box (see Figure \ref{fig:methodology}a) at uniformly spaced Cartesian grid points as the distance from these points to the triangle mesh\footnote{The grid values are used to evaluate, through interpolation, the level-set values and gradients anywhere near or within each block. See \cite{kawamoto_level_2016} for more details on the level-set geometric representation in LS-DEM}.
Figure \ref{fig:configuration}(a) illustrates the triangle mesh representing the surface of the central block with $\theta=10^{\circ}$ and its bounding box in light purple.
\end{comment}

Fig. \ref{fig:TIS}(e) shows $P-\delta$ curves corresponding to four analyses made with the original LS-DEM, using $k_n$=2 GN/mm.
Fig. \ref{fig:TIS}(f) shows $P-\delta$ curves obtained with the most refined surface discretization (Nnodes=68,000) and a set of four evenly-spaced $k_n$'s.
Fig. \ref{fig:TIS}(a) shows that the $P-\delta$ curves scale linearly with both Nnodes and $k_n$, similarly to the two-sphere example shown in Fig. \ref{fig:twoSphereCompression_original}. 
Fig. \ref{fig:TIS}(g) shows three $P-\delta$ obtained with the adapted LS-DEM for varying Nnodes between 17,000-68,000 and $k_n^{*}$=2 GPa/mm.
Fig. \ref{fig:TIS}(h) shows three $P-\delta$ obtained with the adapted LS-DEM, Nnodes=68,000, and a set of four linearly spaced $k_n^{*}$'s.
From the essential overlapping of the three curves in Fig. \ref{fig:TIS}(g) it is clear that the adapted LS-DEM converges with less than Nnodes=17,000.
A fuller illustration of the convergence process is shown and discussed later in relation to Fig. \ref{fig:convergence05}.
From Fig. \ref{fig:TIS}(h), the response in the adapted formulation scales linearly with $k_n^{*}$, in accordance with Eq. (\ref{eq:penetration_stiffness_adapted_LSDEM}).

Comparing Figs. \ref{eq:penetration_stiffness_original_LSDEM}(c,d) and \ref{fig:twoSphereCompression_adapted}(a,b) to Figs. \ref{fig:TIS}(c,d,e,f), respectively, shows that, in both the two-sphere case and the slab example, the dependence of the original and adapted formulations of Nnodes, $k_n$, and $k_n^{*}$ is essentially the same.
This means that neither the divergent nature of the original formulation nor the convergence of the adapted one dependent on the number on grains (blocks) in the problem.

\subsection{The relation between SD and VD refinement} \label{subsection:The relation between SD and VD refinement}

In the context of the original LS-DEM, it is not clear how to determine the relation between the SD refinement and the VD refinement.
The convergence of the adapted LS-DEM provides a clear-cut criterion for this determination, namely that the SD and VD should be refined such that the results do not change with further refinement of either or both together.
Going through such a two-way convergence study could provide a general estimate to how refined the SD should be relative to the VD.
To obtain such an estimate, we consider next the two-way SD-VD convergence in the context of the topologically interlocked slab discussed in section \ref{subsection:A topologically interlocked slab}.

\begin{figure}[H] 
    \includegraphics[width=0.95\textwidth]{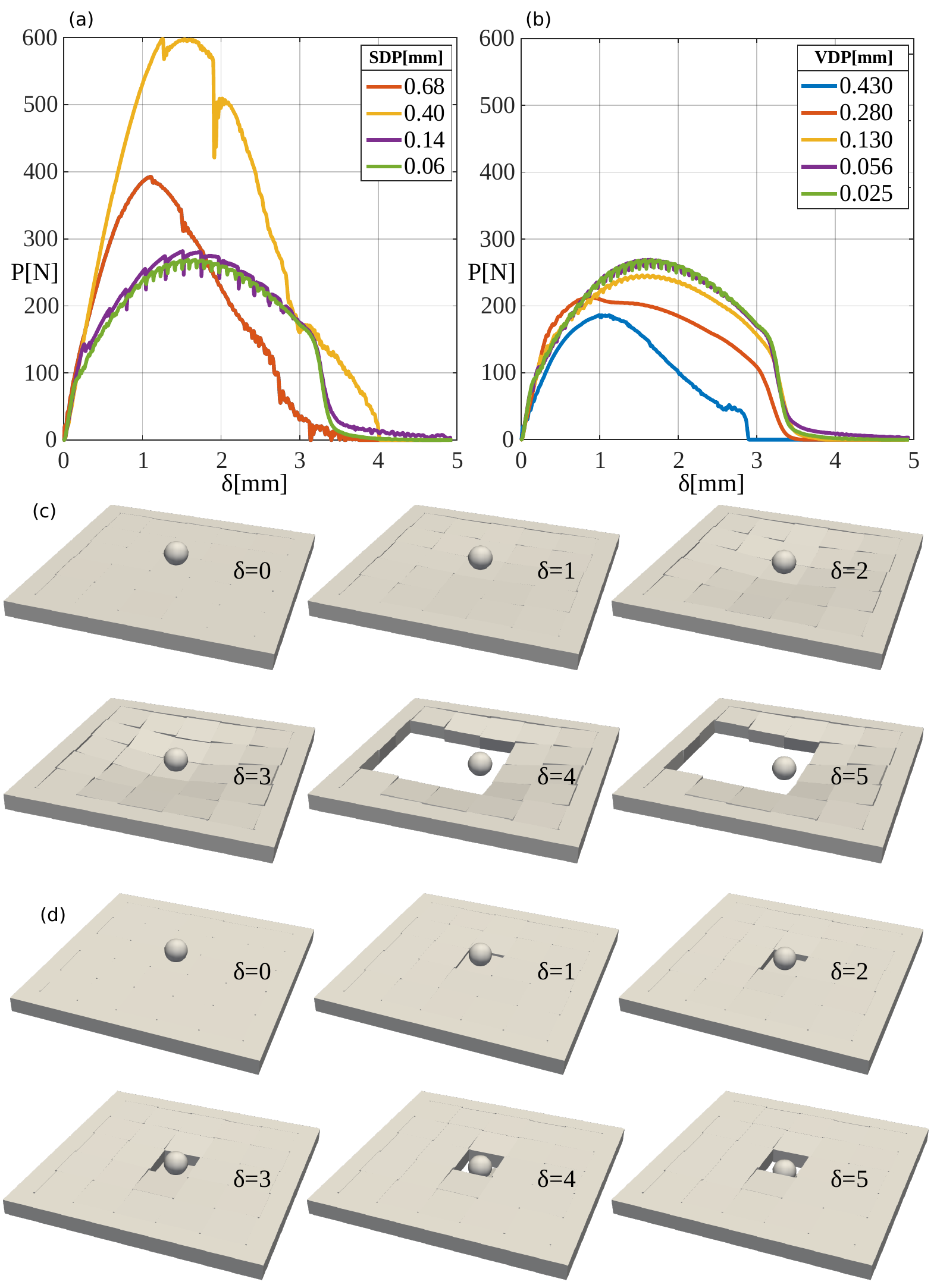} 
    \centering
    \caption{The adapted LS-DEM is fully discretization convergent: $P-\delta$ curve convergence with surface (a) and volume (b) discretization refinement; failure mechanism convergence from coarse (c) to refined (d) discretizations}
    \label{fig:convergence05}
\end{figure}

Fig. \ref{fig:convergence05}(a,b) depict the two-way SD-VD convergence.
Here, we choose the distance between surface nodes as the surface discretization parameter SDP, and the step size in the Level-Set volume grid as the volume discretization parameter VDP.
The reason is that length parameters are a more meaningful measure for quantifying the relation between the SD and VD refinement than the number of nodes, which scales differently for SD and VD with the overall lineal dimensions of the grain.

Fig. \ref{fig:convergence05}(a) shows $P-\delta$ obtained for a series of SDP's at a constant VDP=0.025 mm.
Considering the negligible differences between the two most refined analyses, SDP=0.06 mm (with Nnodes=68,314) was deemed the converged SDP.
We note that the maximal $P$ varies by a factor of approximately 2 as Nnodes grows by more than 2 orders of magnitude, from 420 with SDP=0.68 mm to 68,314 with SDP=0.06 mm.
This speaks to the small sensitivity of the adapted formulation on the SD refinement, which allows to obtain reasonable estimates of the converged behavior with computationally cheaper coarse SDs. 
Fig. \ref{fig:convergence05}(b) shows $P-\delta$ obtained for a series of VDP's at a constant converged SDP=0.06 mm.
Considering the negligible differences between the two most refined analyses, VDP=0.06 mm was deemed the converged SDP.

The fact that the response reached saturation with discretization refinement at SDP=0.06 mm and VDP=0.025 mm allows defining the VDP to SDP ratio $\frac{0.025}{0.06}=0.416$ required for convergence. 
This ratio may serve as a first-order approximation to facilitate future two-way SD-VD convergence studies.

\subsection{The physical significance of convergence} \label{subsection:The physical significance of convergence}

How important is it to use a converged discretization in terms of correctly capturing the mechanical response?
We address this question in the context of the topologically interlocked slabs by comparing the failure mechanism as obtained with the adapted LS-DEM to the experimentally observed one in \cite{mirkhalaf_toughness_2019}

In \cite{mirkhalaf_toughness_2019}, the failure was governed by sliding of the central block and its eventual falling off.
The failure mechanisms as obtained with the adapted LS-DEM are depicted in Figs. \ref{fig:convergence05}(c,d), through six snapshots corresponding to $\delta=0,1,...,5$ mm.
Fig. \ref{fig:convergence05}(c) shows that with the coarse discretization (SDP=0.68 mm and VDP=0.43 mm), the failure mechanism is qualitatively different from the experimentally observed one. 
The blocks stick to one another throughout most of the response and several of the central blocks end up falling-off.
In contrast, Fig. \ref{fig:convergence05}(d) shows that the converged discretization correctly captures the experimentally observed failure.
This speaks to the importance of having a converged LS-DEM formulation and working with converged discretizations for capturing the correct mechanical response. 

%\subsection{The relative importance of the LSDEM adaptations in different types of problems} \label{subsection:The relative importance of the LSDEM adaptations in different types of problems}

\section{Summary and outlook} \label{sec:Summary and outlook}
In this paper, we have shown that the original LS-DEM diverges upon surface discretization refinement, even in a simple case of central compression between two spheres.
We have further shown that this divergence is due to the force-based formulation of nodal contact interactions and the commensurate linear scaling of the grain penetration stiffness with the number of contacting surface nodes.
Based on this, we have replaced the original contact formulation with a continuum-based one wherein the nodal contact interactions are traction-based.
We have shown that the adapted LS-DEM is fully discretization convergent in different kinds of applications including the novel application of LS-DEM to structural analysis.
Based on full two-way convergence study for the surface and the volume discretizations, we have derived the volume-to-surface degree of refinement necessary for convergence, and we have shown the discretization convergence is physically significant in that it is necessary to correctly capture the mechanical response.
The convergent nature of the adapted formulation paves the way to broaden the scope of LS-DEM applications to structural analysis of structures made of discrete building blocks.
In such applications, the divergence of the original LS-DEM is most pronounced and therefore the convergence of the adapted LS-DEM is essential to obtain quantitatively meaningful results.
Lastly, the fact that the additional computational cost of the adapted formulation is effectively zero, warrants its general use in future LS-DEM applications.
\printbibliography

\end{document}